\documentclass{article}

\usepackage[latin1]{inputenc}
\usepackage{graphicx}
\usepackage{epsfig}
\newtheorem{theorem}{Theorem}[section]

\newtheorem{remark}{Remark}[section]
\newtheorem{definition}{Definition}[section]

\newfont{\bb}{msbm10 at 10pt}
\def\r{\hbox{\bb R}}

\def\n{\hbox{\bf n}}

\begin{document}

\title{A criterion on instability of rotating cylindrical surfaces}
\author{Rafael L\'opez\footnote{Partially
supported by MEC-FEDER
 grant no. MTM2007-61775.}\\
Departamento de Geometr\'{\i}a y Topolog\'{\i}a\\
Universidad de Granada\\
18071 Granada. Spain\\
rcamino@ugr.es}
%\subjclass{Primary 53C42; Secondary 53C20}
%\keywords{Rotating surface, stability,  Rayleigh condition}
\date{}

\maketitle
\begin{abstract} We consider a column of a rotating stationary surface  in Euclidean space. We obtain  a value $l_0>0$ in such way that if the length $l$ of column satisfies $l>l_0$, then the surface is instable. This extends, in some sense, previous results due to Plateau and Rayleigh for  columns of surfaces with constant mean curvature.
\end{abstract}
%%%%%%%%%%%%%%%%%%%%%%%%%%%%%%%%%%%%%%%%%%%%%%%%%%%
\section{Statement of the problem and the result}
%%%%%%%%%%%%%%%%%%%%%%%%%%%%%%%%%%%%%%%%%%%%%%%%

Consider a column of a circular cylinder of radius $r$.   One asks what  lengths $l$ of columns determine stable cylinders (as surfaces with constant mean curvature). Rayleigh \cite{ra} proved that the column is stable iff $l<2\pi r$ (see also \cite{bc}). This was experimentally  observed by Plateau in 1873. In this work, we give  a similar result for rotating cylindrical surfaces. Firstly, we present the physical framework necessary for the statements of our results.

We consider the steady rigid rotation of
an homogeneous incompressible fluid drop in absence of gravity, which is surrounded by a
rigidly rotating incompressible fluid. In mechanical
equilibrium, we say that the drop is a rotating stationary drop, or simply, a rotating drop. Rotating  drops have been the
subject of intense study beginning from the work of  Plateau \cite{pl}, because they could be models in other areas of physics, such as,
astrophysics, nuclear physics, or fluid dynamics. In the case of our interest, the surfaces are of cylindrical type, that is, ruled surfaces whose rulings are parallel to the axis of rotation. These surfaces have appeared in a number of physical  settings, usually referred in the literature as the problem of the shape of a \emph{Hele-Shaw cell}. For example, if the fluid rotates with respect to a vertical line, we consider that the liquid lies confined in a narrow slab determined by two parallel horizontal plates. Then we can suppose that the shape of the drop is invariant by vertical displacements. We refer the reader to  \cite{cmco,cso,sc} and references therein for further details  and  a more precise statements of the physical problem.

From a mathematical viewpoint,  we  assume the liquid drop rotates about the
$x_3$-axis with a  constant angular velocity $\omega$, where $(x_1,x_2,x_3)$ are the usual coordinates of Euclidean space $\r^3$.
Let $\rho$ be the constant density of the fluid of the drop and denote by  $W$  the bounded open set
in $\r^3$, which is the region occupied by the rotating drop. We set $S=\partial W$ as the free interface between the drop and the ambient liquid  and that we suppose to be a smooth boundary surface. The energy of this mechanical system is given by
$$E=\tau|S|-\frac12\rho\omega^2\int_{W} r^2 dx,$$
where $\tau$ stands for the surface tension on $S$, $|S|$ is the surface area of $S$ and $r=r(x)=\sqrt{x_1^2+x_2^2}$ is the distance
from a point $x$ to the $x_3$-axis.
The term  $\tau|S|$ is the surface energy of the drop and
$\frac12\rho\omega^2\int_W r^2 dx$ is the potential energy associated
with the centrifugal force. We assume that the volume $V$ of the drop remains constant while  rotates. The equilibrium shapes of such a drop are governed by the Laplace equation
$$2\tau H(x)=-\frac12\rho\omega^2 r^2+\lambda\hspace*{1cm}\mbox{$x\in S$},$$
where $\lambda$ is a constant depending on the volume constraint.
As consequence, the mean curvature of the interface $S$ satisfies an equation of type
$$H(x_1,x_2,x_3)=a r^2+b,$$
where $a,b\in\r$. We say then that $S$ is a  \emph{rotating stationary surface}. When the liquid does not rotate ($a=0$, or $\omega=0$ in the Laplace equation),  the interface is modeled by a surface with constant mean curvature.

Consider $M$ a cylindrical surface whose generating curve $C$  is a closed curve that lies in the $x_1x_2$-plane. The rulings of the surface are straight-lines parallel to the $x_3$-axis and we parametrize $M$ as
$x(t,s)=\alpha(s)+t(0,0,1)$, $t,s\in\r$, where $\alpha$ is a parametrization of $C$. We always orient $M$ so the Gauss map points inside (this means that $\alpha$ runs in the counterclockwise direction). We consider \emph{columns} of cylindrical surfaces, that is,   the range of the parameter $t$ in the above parametrization of $M$ is an interval of length $l$.

It is important to point out that there exist cylindrical rotating surfaces whose generating curve is \emph{closed}. Exactly, the variety of shapes of closed curves are as follows. We begin with a circle corres\-ponding with zero angular velocity. As we increase the angular velocity,
the curve changes through a sequence of shapes which evolved from non-self-intersecting curves with fingered morphology and high symmetry for slow rotation to self-intersecting curves with a richness of  symmetric properties for fast angular velocity. A complete description of these curves can seen in \cite{lom}.

Recall that a rotating stationary surface is called \emph{stable} if   the second variation of the energy of the surface is non-negative for every smooth deformation of the surface which fixes the enclosed volume and fixes the boundary.  In this paper we will prove the following

\begin{theorem} \label{main} Let $M$ be a cylindrical rotating surface of column length $l$ that satisfies the Laplace equation $2H(x)=ar^2+b$. Let  $C$ be the planar generating curve of $M$  in the $x_1x_2$-plane. Let  $\Omega$ be the bounded domain by $C$ and denote by $L$ and $|\Omega|$ the length of $C$ and the area of $\Omega$ respectively. If $a>0$ and  $M$ is stable, then
$$l\leq \frac{\pi}{2}\sqrt{\frac{L}{a|\Omega|}}.$$
\end{theorem}

In relation to this,   stability of rotating stationary  surfaces   has been studied for rotationally symmetric configurations (invariant by rotations with respect to the $x_3$-axis). See   \cite{agu,bs,bs2,ca,mc,ws}.

%%%%%%%%%%%%%%%%%%%%%%%%%%%%%%%%%%%%%%%%%%%%%%%%%%%
\section{Stability of rotating surfaces}
%%%%%%%%%%%%%%%%%%%%%%%%%%%%%%%%%%%%%%%%%%%%%%%%

In this section, we give some preliminaries and definitions, and we fix the notation.
Let $x:M\rightarrow\r^3$ be an immersed oriented compact surface and consider a smooth variation $X$ of $x$, that is, a smooth map $X:(-\epsilon,\epsilon)\times M:M\rightarrow\r^3$  such that, by setting $x_t=X(t,-)$, we have $x_0=x$ and $x_t-x\in C^{\infty}_0(M)$. We require that the volume $V(t)$ of each immersion $x_t$  remains constant throughout the variation. Let $u\in C^{\infty}(M)$ be the normal component of the variational vector field of $x_t$,
$$u=\langle\frac{\partial x_t}{\partial t}{\bigg|}_{t=0},N\rangle,$$
where $N$ is the Gauss map of $x$.
We take $E(t):=E(x_t)$ the value of the energy for each immersion $x_t$ and $V(t)$ the corresponding volume. The first variation of $E$ and $V$ at $t=0$ is given by
$$E'(0)=-\int_M\bigg(\frac12\rho\omega^2 r^2+2\tau H\bigg) u\ dM\hspace*{1cm}\mbox{and}\hspace*{1cm}V'(0)=\int_M u\ dM,$$
respectively. Here $H$ is the mean curvature of the immersion $x$.
By the method of Lagrange  multipliers,  the first variation of $E$ at $t=0$ is to
be zero relative for all volume preserving variations if there is a constant $\lambda$ so that
$E'(0)+\lambda V'(0)=0$.  This yields the
condition
$$2\tau H=-\frac12\rho\omega^2 r^2+\lambda\hspace*{1cm}\mbox{on $M$}.$$
See \cite{we2} for details.  Thus the mean curvature $H$ satisfies an equation of type
\begin{equation}\label{laplace}
2H(x)=a r^2+b,\hspace*{1cm}a,b\in\r.
\end{equation}
Denote by $\{E_1,E_2,E_3\}$ the canonical orthonormal basis  of  $\r^3$:
$$E_1=(1,0,0),\hspace*{.5cm}E_2=(0,1,0),\hspace*{.5cm}E_3=(0,0,1).$$
We write $N_i=\langle N,E_i\rangle$ and $x_i=\langle x,E_i\rangle$
for $1\leq i\leq 3$. We  obtain  the second variation of the energy  in order to give the notion of stability. In this section, we present a derivation of the second variation of the energy for a rotating surface that follows ideas of M. Koiso and B. Palmer \cite{kp} (a general formula of the second variation was obtained by Wente \cite{we}; see also \cite{mc}). Assume that $x:M\rightarrow\r^3$ is a critical point of $E$ and we write the second variation of the functional $E$ in the form
$$E''(0)=-\int_M u\cdot L[u]\ dM,$$
where $L$ is a linear differential operator acting on  the normal component, which we want to find.
Since the translations in the $E_3$  are symmetries of the energy functional $E$, then $L[N_3]=0$. The same occurs for   the rotation with respect to the $x_3$-axis,  $L[\psi]=0$,  where  $\psi=\langle x\wedge N,E_3\rangle$ and $\wedge$ is the vector cross product of $\r^3$. We now compute  $L[N_i]$ and $L[\psi]$. The tension field of the Gauss map satisfies
$$\Delta N+|\sigma|^2 N=-2\nabla H,$$
where $\Delta$ is the Laplacian in the metric induced by $x$,
$\sigma$ is the second fundamental of the immersion  and $\nabla$ is
the covariant differentiation. Since $2H=ar^2+b$, we have
\begin{equation}\label{dh}
2\nabla H=2a(x_1\nabla x_1+x_2\nabla x_2)=2a\bigg( x_1 E_1+x_2 E_2-(h-x_3 N_3)N\bigg).
\end{equation}
Here $h=\langle N,x\rangle$ stands for the support function of $M$.
In particular,
\begin{equation}\label{ene2}
\Delta N_3+\bigg(|\sigma|^2-2a(h-x_3 N_3)\bigg)N_3=0
\end{equation}
We now take the function $\psi$. In general, we have
$$\Delta \psi+|\sigma|^2 \psi=-2\langle \nabla H,E_3\wedge x\rangle.$$
It follows from (\ref{dh}) that
$$-2\langle\nabla H,E_3\wedge x\rangle=2a(h-x_3 N_3)\langle N,E_3\wedge x\rangle=2a(h-x_3 N_3)\psi.$$
Therefore,
\begin{equation}\label{ene3}
\Delta \psi+\bigg(|\sigma|^2-2a(h-x_3 N_3)\bigg)\psi=0.
\end{equation}
As a consequence of (\ref{ene2}) and (\ref{ene3}),
\begin{equation}\label{ele}
L=\Delta+|\sigma|^2-2a(h-x_3 N_3).
\end{equation}
\begin{definition} Consider $x:M\rightarrow\r^3$ be a smooth immersion that satisfies the Laplace equation (\ref{laplace}).
We say that $x$ is stable if
\begin{equation}\label{second}
\Phi(u):=-\int_M u\bigg(\Delta u+\big(|\sigma|^2-2a(h-x_3 N_3)\big)u\bigg)\ dM\geq 0
\end{equation}
for all $u\in C_0^{\infty}(M)$ such that
$$\int_M u\ dM=0.$$
\end{definition}

%%%%%%%%%%%%%%%%%%%%%%%%%%%%%%%%%%%%%%%%%%%%%%%%%%%
\section{Proof of the Theorem}
%%%%%%%%%%%%%%%%%%%%%%%%%%%%%%%%%%%%%%%%%%%%%%%%

Consider now that $M$ is a column of a cylindrical rotating surface of length $l$ whose generating curve $\alpha=\alpha(s)$  is a  closed curve in the $x_1x_2$-plane with $s$ the arc-length parameter. We denote by $\n(s)$ the unit normal vector to $\alpha$, that is, $\n(s)=J\alpha'(s)$, where $J$ is the $\pi/2$-rotation of the $x_1 x_2$-plane along the counterclockwise direction. In particular, $\{\alpha'(s),\n(s),E_3\}$ is an oriented orthonormal basis of $\r^3$. Moreover,  $\alpha''(s)=\kappa(s)\n(s)$, where $\kappa$ is the curvature of $\alpha$.

We parametrize $M$  as $x(t,s)=\alpha(s)+t E_3$, $t\in [0,l]$ and $s\in [0,L]$, where $L$ denotes  the length of  $\alpha$. The first fundamental form of $M$ is $I=dt^2+ds^2$ and the Gauss map is
$N(t,s)=E_3\wedge \alpha'(s) $.  A computation of the mean curvature $H$ gives
$$2H(t,s)= \langle E_3\wedge \alpha'(s),\alpha''(s)\rangle= \kappa(s).$$
We calculate each one of the terms of (\ref{ele}). Since the first fundamental form of $M$ is the Euclidean metric,  $\Delta$ coincides with the Euclidean Laplacian operator $\Delta_0=\partial^2_{tt}+\partial^{2}_{ss}$. The term $|\sigma|^2$ is $\kappa(s)^2$ and
$$h-N_3x_3=\langle N(t,s),x(t,s)\rangle=-\langle\alpha(s)\wedge \alpha'(s),E_3\rangle.$$
Consider smooth functions on $M$ written in the form $u(t,s)=f(t)g(s)$. Then $u$ is a test-function for $\Phi$ if $f(0)=f(l)=0$, $g(0)=g(L)=g'(0)=g'(L)$ and
$$\int_0^l f(t)\ dt=0\hspace*{1cm}\mbox{or}\hspace*{1cm}\int_0^L g(s)\ ds=0.$$
From the above computations, $\Phi(u)$ in (\ref{second}) writes as
$$\Phi(u)=-\int_0^L\int_0^l f(t)g(s)\bigg( f''(t)g(s)+f(t)g''(s)+\big(\kappa(s)^2+2a\langle\alpha(s)\wedge \alpha'(s),E_3\rangle\big) f(t)g(s)\bigg)\ dt\ ds.$$

Let us take $f(t)=\sin{\big(\frac{\pi t}{l}\big)}$ and $g_i(s)=\langle \alpha'(s),E_i\rangle$, $i=1,2$. The functions $g_i$ correspond, up a sign, with the functions $N_i$, $i=1,2$. Then $\int_0^L g_i(s)\ ds=0$. Let $u_i(t,s)=f(t)g_i(s)$, $i=1,2$. If $M$ is stable, then
$\Phi(u_i)\geq 0$, $i=1,2$. Let do the computations of $\Phi(u_i)$. As $f''(t)=-\frac{\pi^2 }{l^2}f(t)$, we have that $\Phi(u_i)\geq 0$ iff
\begin{equation}\label{esta1}
\frac{\pi^2 }{l^2}\int_0^L  g_i(s)^2\ ds -\int_0^L g_i(s)g_i''(s)\ ds\geq \int_0^L \bigg(\kappa(s)^2+2a\langle\alpha(s)\wedge \alpha'(s),E_3\rangle\bigg) g_i(s)^2 \ ds.
\end{equation}
An integration by parts leads to $\int_0^Lg_i(s)g_i''(s)\ ds=-\int_0^L g_i'(s)^2\ ds$. We note that $g_i'(s)=\kappa(s)\langle\n(s), E_i\rangle$. Therefore,
$$g_1(s)^2+g_2(s)^2=1,\hspace*{1cm}g_1'(s)^2+g_2'(s)^2=\kappa(s)^2.$$
 Summing the two expressions of (\ref{esta1}) for $i=1,2$, we have
$$\frac{\pi^2}{l^2}L\geq 2a \int_0^L\langle\alpha(s)\wedge\alpha'(s),E_3\rangle\ ds.$$
As the curve $\alpha$ runs along the  counterclockwise direction, then
$$\int_0^L\langle\alpha'(s)\wedge\alpha'(s),E_3\rangle\ ds=2|\Omega|,$$
where $\Omega$ is the bounded domain determined by $\alpha$ in the $x_1x_2$-plane.
Therefore
$$l^2\leq \frac{\pi^2 L}{4a |\Omega|}$$
which yields the result.

One can check that the function $u(t,s)=\sin(\frac{2\pi t}{l})$ is a test-function for the functional $\Phi$ since
$\int_0^l f(t)\ dt=0$ and $f(0)=f(l)=0$. Doing similar computations, we obtain that if $M$ is stable then
$$\frac{4\pi^2}{l^2}L\geq \int_0^L\kappa(s)^2\ ds+4a|\Omega|=\int_0^L (a r^2+b)^2\ ds+ 4a|\Omega|> 4a|\Omega|.$$
In particular, if $a>0$, $l^2<\frac{\pi^2}{a|\Omega|}L$. This inequality is weaker the one obtained in Theorem \ref{main}. However, if $b\geq 0$, then one can follow estimating with $\int_0^L (ar^2+b)^2ds\geq L b^2$. Then
$$l^2\leq\frac{4\pi^2 L}{L b^2 +4a|\Omega|}.$$

\begin{theorem}\label{main2}Under the same hypothesis and notation as in Theorem \ref{main}, if  $a\geq 0$, $b>0$ and  $M$ is stable, then
$$l\leq 2\pi\sqrt{\frac{L}{4a|\Omega|+L b^2 }}.$$
\end{theorem}

This result improves the one obtained in Theorem \ref{main} if $\frac{a}{b^2}<\frac{L}{12|\Omega|}$.

\begin{remark} Consider $M$ a circular cylinder of radius $r$. Then $M$ can be viewed in a twofold sense. First, as a  surface with constant mean curvature $H=\frac{1}{2r}$ and second, as a rotating surface, with $a=1/r^3$ and $b=0$. Thus we have two different notions of stability for the same surface.  As it was mentioned at the beginning of this work, it is known that  $M$ is stable as a constant mean curvature surface iff $l< 2\pi r$. If we compare with our result, Theorem \ref{main} says that if $M$ is stable as a rotating surface then
$l\leq \frac{\pi r}{\sqrt{2}}$. As a consequence, there exist columns of circular cylinders that are instable as   rotating surfaces but they are stable as surfaces with constant mean curvature.
\end{remark}

\begin{remark}  Consider $M$ a circular cylinder of radius $r$ and length $l$ again. According with our choice of the orientation of the generating curve,  $M$ satisfies the Laplace equation (\ref{laplace}) for $a=0$ and $b=1/r$. The estimate of the length $l$ for stable columns of cylinders given in Theorem \ref{main2}  reads now as $l\leq 2\pi r$, re-discovering the value obtained by Rayleigh.
\end{remark}

%%%%%%%%%%%%%%%%%%%%%%%%%%%%%%%%%%%%%%%%%%

\end{document}